\documentclass[]{amsart}
\usepackage[foot]{amsaddr}
\usepackage[a4paper,top=2.5cm,bottom=3cm,inner=3cm,outer=3cm]{geometry}

\usepackage[utf8]{inputenc}
\usepackage{tikz}

\usepackage{amsmath,amsthm}
\usepackage{amsfonts,amssymb,mathrsfs}

\usepackage[bookmarks=false,breaklinks=true]{hyperref}
\usepackage{float}
\usepackage{graphicx}
\usepackage{color}

\usepackage{epigraph}
\definecolor{myurlcolor}{rgb}{0,0,0.7}
\usepackage{hyperref}
\hypersetup{colorlinks,
linkcolor=myurlcolor,
citecolor=myurlcolor,
urlcolor=myurlcolor}
\usepackage[mathscr]{euscript}

\input xy
\xyoption{all}

\renewcommand{\H}{\mathcal{H}}  
\newcommand{\h}{\mathfrak{h}}   

\newcommand{\tr}{\mathrm{tr}}    

\newcommand{\Q}{{\mathbb Q}}  
\newcommand{\R}{{\mathbb R}}  
\newcommand{\C}{{\mathbb C}}  
\newcommand{\E}{{\mathbb E}}  

\newcommand{\x}{\mathbf{x}}

\newcommand{\define}[1]{{\bf \boldmath{#1}}}

\theoremstyle{definition}

        \newcommand{\be}{\begin{equation}}
        \newcommand{\ee}{\end{equation}}
        \newcommand{\ba}{\begin{eqnarray}}
        \newcommand{\ea}{\end{eqnarray}}
        \newcommand{\ban}{\begin{eqnarray*}}
        \newcommand{\ean}{\end{eqnarray*}}
        \newcommand{\barr}{\begin{array}}
        \newcommand{\earr}{\end{array}}

\begin{document}
\title{The Hexagonal Tiling Honeycomb}
\author[Baez]{John C.\ Baez} 
\address{Department of Mathematics, University of California, Riverside CA, 92521, USA}
\date{November 2024}
\maketitle

\begin{center}
\frame{\includegraphics[width = 35 em]{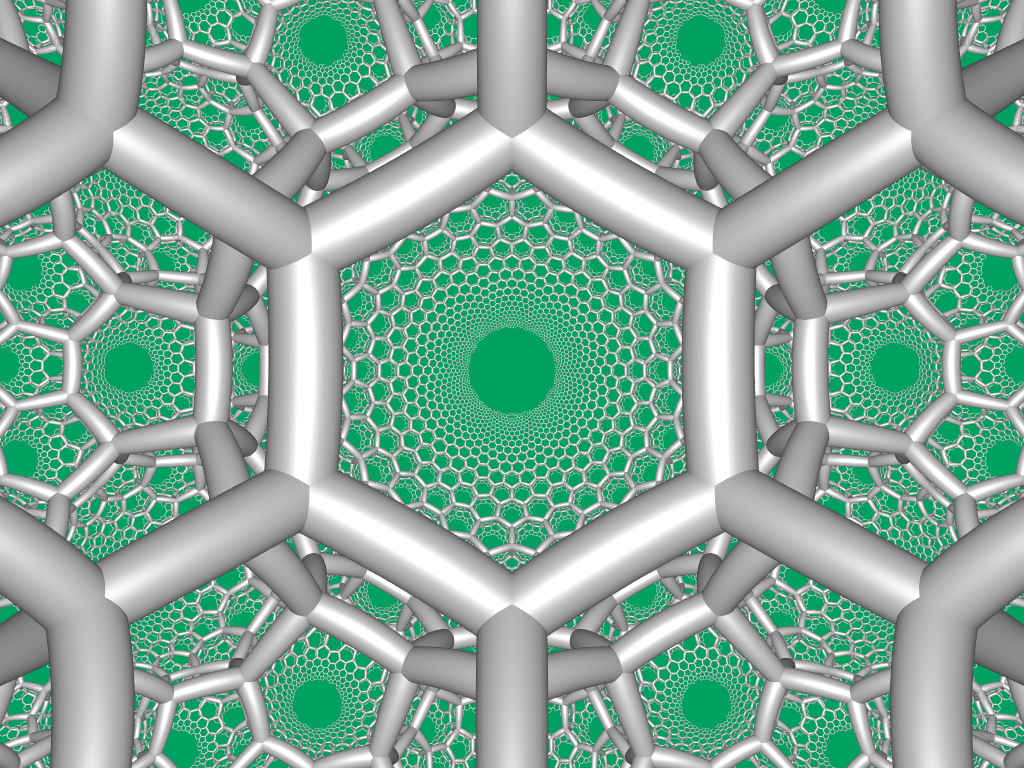}}
\end{center}

\vskip 1em

This picture by Roice Nelson shows the `hexagonal tiling honeycomb'.  But what is it?  Roughly speaking, a honeycomb is a way of filling 3d space with polyhedra.  The most familiar one is the usual way of filling Euclidean space with cubes.  This cubic honeycomb is denoted by the symbol $\{4,3,4\}$, because a square has 4 edges, 3 squares meet at each corner of a cube, and 4 cubes meet along each edge of this honeycomb.  We can also define honeycombs in hyperbolic space, which is a 3-dimensional Riemannian manifold with constant negative curvature.   For example, there is a hyperbolic honeycomb denoted $\{4,3,5\}$, which is a way of filling hyperbolic space with cubes where 5 meet along each edge.

Coxeter has classified the most symmetrical hyperbolic honeycombs \cite{Coxeter}, called the `regular' ones, and the hexagonal tiling honeycomb is one of these.  But it does not contain polyhedra of the usual sort, with finitely many faces.  Instead, it contains flat Euclidean planes embedded in hyperbolic space, each of which is tiled by regular hexagons.   These tiled planes can be seen as generalized polyhedra with \emph{infinitely many} faces.   The symbol for the hexagonal tiling honeycomb is $\{6,3,3\}$, because a hexagon has 6 edges, 3 hexagons meet at each corner in a plane tiled by regular hexagons, and 3 such planes meet along each edge of this honeycomb.

The hexagonal tiling honeycomb shows up naturally if we try to discretize spacetime while preserving as much symmetry as we can.   In special relativity, \define{Minkowski spacetime} is $\R^4$ equipped with the nondegenerate bilinear form 
\[   (t,x,y,z) \cdot (t',x',y',z') = tt' - xx' - yy' - zz',  \] 
usually called the \define{Minkowski metric}.
Hyperbolic space sits inside Minkowski spacetime as the hyperboloid of points $\x = (t,x,y,z)$ with $\x \cdot \x = 1$ and $t > 0$.   Equivalently, we can think of Minkowski 
spacetime as the space $\h_2(\C)$ of $2 \times 2$ hermitian complex matrices, using the fact that every such matrix is of the form 
\[    A =  \left( \begin{array}{cc} t + z & x - i y \\ x + iy & t - z \end{array} \right)\! \]
and $\det(A) =  t^2 - x^2 - y^2 - z^2$.    In these terms, hyperbolic space is the hyperboloid
\[        \H = \left\{A \in \h_2(\C) \, \vert \, \det A = 1, \,  \tr(A) > 0 \right\}   .\]

How can we construct the hexagonal tiling honeycomb inside $\H$?   Sitting in the complex numbers we have the ring $\E$ of \define{Eisenstein integers}: complex numbers of the form $a + b \omega$ where $a,b$ are integers and $\omega$ is a nontrivial cube root of $1$, say $\omega = e^{2 \pi i / 3}$.    This lets us define a lattice in Minkowski spacetime, say $\h_2(\E)$, consisting of $2 \times 2$ hermitian matrices with entries that are Eisenstein integers.  This lattice can be seen as a discretized version of Minkowski spacetime.  Then comes a minor miracle: the points at the centers of hexagons in the hexagonal tiling honeycomb are precisely those points in the lattice $\h_2(\E)$ that lie on the hyperboloid $\H$.  For two proofs see \cite{B}.

The hexagonal tiling honeycomb also arises in algebraic geometry.   An \define{abelian variety} is a complex projective variety that is also an abelian group in a compatible way.  The most famous are the 1-dimensional ones,  called \define{elliptic curves}.  Any elliptic curve is of the form $\C/L$ for some lattice $L \subset \C$.  We can get a highly symmetrical elliptic curve from the Eisenstein integers by forming the quotient $\C/\E$.   We can then form a 2-dimensional abelian variety by taking the product $\C/\E \times \C/\E \cong \C^2/\E^2$.   

The set of isomorphism classes of complex line bundles over any complex projective variety becomes an abelian group thanks to our ability to tensor line bundles.   This group, called the \define{Picard group}, has a natural topology, and it typically has many connected components.  The set of connected components is an abelian group in its own right, called the \define{N\'eron--Severi group}.   You can think of this as the group of equivalence classes of line bundles where two count as equivalent if one can be deformed to the other.

Thanks to some beautiful theorems on abelian varieties \cite[Chap.\ 5]{BL}, the N\'eron--Severi group of $\C^2/\E^2$ is isomorphic to $\h_2(\E)$.   The points $A \in \h_2(\E)$ with $\tr(A) > 0$ and $\det(A) > 0$ come from `ample' line bundles, which are roughly those having enough sections to provide an embedding of the variety into projective space.    Among these points, those with $\det(A) = 1$ correspond to so-called `principal polarizations', which play an important role in the study of abelian varieties.   But these points are precisely the centers of the hexagons in the hexagonal tiling honeycomb!   

Thus, the picture above gives a vivid example of some concepts from algebraic geometry.   But it is also part of a larger story relating algebraic integers in the fields $\Q[\sqrt{-n}]$ to regular honeycombs \cite{Garner} and their symmetry groups, called `Coxeter groups' \cite{JW,MSW}.

\end{document}